\documentclass[graybox]{svmult}

\usepackage{mathptmx}       % selects Times Roman as basic font
\usepackage{helvet}         % selects Helvetica as sans-serif font
\usepackage{courier}        % selects Courier as typewriter font
\usepackage{type1cm}        % activate if the above 3 fonts are
\usepackage{makeidx}         % allows index generation
\usepackage{graphicx}        % standard LaTeX graphics tool
\usepackage{multicol}        % used for the two-column index
\usepackage[bottom]{footmisc}% places footnotes at page bottom

\newcommand{\reals}{\mbox{${\textrm I}\!{\textrm R}$}}

\usepackage{amsfonts}
\makeindex             % used for the subject index
\begin{document}

\title*{Negative Cycle Separation in Wireless Network Design%
\thanks{This work was partially supported by the \emph{German Federal Ministry of Education and Research} (BMBF), project \emph{ROBUKOM}, grant 03MS616E, and by the \emph{Italian Ministry of Education, University and Research} (MIUR), project \emph{``Design, optimization and analysis of next-generation converged networks''}, grant PRIN 2008LLYXFS.
\newline
\newline
This is the authors' final version of the paper published in Pahl J., Reiners T., Voß S. (eds) Network Optimization. Lecture Notes in Computer Science, vol 6701, pp. 51-56. Springer, Berlin, Heidelberg. DOI: 10.1007/978-3-642-21527-8\_7
\newline
The final publication is available at Springer via http://dx.doi.org/10.1007/978-3-642-21527-8\_7}%
}
\author{Fabio D'Andreagiovanni,
Carlo Mannino and Antonio Sassano}
\institute{Fabio D'Andreagiovanni \at Konrad-Zuse-Zentrum f\"ur Informationstechnik Berlin (ZIB),
Takustr. 7, D-14195 Berlin, Germany,
\email{d.andreagiovanni@zib.de}
\and Carlo Mannino \at Department of Computer and System Sciences, Sapienza Universit\`a di Roma, via Ariosto 25, 00185 Roma, Italy, \email{mannino@dis.uniroma1.it}
\and Antonio Sassano \at Department of Computer and System Sciences, Sapienza Universit\`a di Roma, via Ariosto 25, 00185 Roma, Italy,
\email{sassano@dis.uniroma1.it}}
\maketitle

\abstract{The Wireless Network Design Problem (WND) consists in choosing values of radio-electrical parameters of transmitters of a
wireless network, to maximize network coverage.
We present a pure 0-1 Linear Programming formulation for the WND that may contain an exponential number of constraints. Violated inequalities of this formulation are hard to separate both theoretically and in practice. However, a relevant subset of such inequalities can be separated more efficiently in practice and can be used to strengthen classical MILP formulations for the WND. Preliminary computational experience confirms the effectiveness of our new technique both in terms of quality of solutions found and provided bounds.}

\section{Introduction}
\label{sec:intro}

Wireless networks have shown a rapid growth over the past two
decades and now play a key role in new generation
telecommunications networks. Scarce radio resources, such as
frequencies, have rapidly became congested and the need for more
effective design methods arose. A general planning problem
consists in establishing the radio-electrical parameters (e.g., power emission and frequency) of the
transmitters of a wireless network so as to maximize the overall
network coverage. To present our original contribution, in this paper we focus only on establishing power emissions. This is actually a basic problem in all wireless planning contexts that can be easily extended by introducing additional elements, such as frequencies \cite{DA10, DAMaSa10}.

For our purposes, a wireless network can be described as a set of
transmitters $B$ distributing a telecommunication service to a set
of receivers $T$. Each transmitter $b\in B$ emits a radio signal with power $p_b\in [0,P_{\max}]$. The power $p_b(t)$ that receiver $t$ gets from transmitter $b$ is proportional to the emitted power $p_b$ by a
factor $a_{tb} \in [0,1]$, i.e. $p_b(t) = a_{tb} \cdot p_b$, commonly called \emph{fading coefficient}. Among the signals received from transmitters in $B$, receiver $t$ can select a {\em reference signal} (or {\em server}), which is the one carrying the service. All the other signals are interfering.

A receiver $t$ is regarded as served by the network, specifically
by server $\beta \in B$, if the ratio of the serving power
to the sum of the interfering powers ({\em signal-to-interference
ratio} or \emph{SIR}) is above a threshold $\delta$ \cite{Ra01},
(\emph{SIR threshold}), whose value depends on the technology
and the desired quality of service:
\begin{equation}
\label{eq:firstSIRinequality} \frac{a_{t \beta} \cdot
p_{\beta}}{\mu + \sum_{b \in B \setminus\{\beta\}} a_{tb}
\cdot p_b} \geq \delta \qquad
\end{equation}

\noindent where the system noise $\mu > 0$   is assimilated to an
interfering signal with fixed (very low) power emission.

For every $t\in T$, we have one inequality of type
(\ref{eq:firstSIRinequality}) for each potential server $\beta\in
B$: in particular, we denote by $SIR(t,b)$ the inequality
(\ref{eq:firstSIRinequality}) associated with receiver $t$ and
server $b$. Receiver $t$ is served if at least one of these
inequalities is satisfied or, equivalently, if the following
disjunctive constraint is satisfied:
\begin{equation}\label{eq:disjunctive-SIR}
\bigvee_{\beta \in B} \left( a_{t \beta} \cdot p_{\beta} - \delta \cdot \sum_{b \in B\setminus\{\beta\}}
a_{tb} \cdot p_b \geq \delta \cdot \mu \right)
\end{equation}

\noindent
Each linear inequality of the above disjunction
is obtained by simple algebra from the SIR
expression (\ref{eq:firstSIRinequality}).

If each receiver $t \in T$ is associated to a value $r_t > 0$
that expresses revenue obtained by serving $t$, the
\emph{Wireless Network Design Problem (WND)} consists in
setting the power emission of each transmitter $b \in B$ and
the server of each receiver in $t \in T$ with the aim of maximizing the overall revenue of served receivers.

\section{A pure 0-1 Linear Programming formulation for the WND}

The WND is often approached by solving a suitable Mixed-Integer
Linear Program (MILP): first, a binary variable $x_{tb}$ is
introduced for every $t\in T$, $b \in B$, with $x_{tb} = 1$ if and
only if $b$ serves $t$; then, variables $x_{tb}$ are used to
replace each disjunction (\ref{eq:disjunctive-SIR}) with a set of
$|B|$ linear constraint, that, however, include large positive
constants, the notorious \emph{big-M coefficients}
\cite{DAMaSa10,KeOlRa10}. The (linear) objective function aims to maximize
the overall revenue from coverage, i.e. $\max \sum_{t\in T}
\sum_{b \in B} r_t \cdot x_{tb}$ and requires the additional
constraints:
\begin{equation}\label{eq:one-server}
\sum_{b \in B} x_{tb} \leq 1 \qquad t\in T
\end{equation}

\noindent
to ensure that each receiver is associated to at most one server. A vector $x\in \{0,1\}^{T\times B}$ satisfying
(\ref{eq:one-server}) is a {\em server assignment}.

The resulting MILP presents severe drawbacks,  highlighted in
several works, e.g. \cite{DAMaSa10,KeOlRa10,MaMaSa09}. First, the
coefficients in the SIR inequalities may vary over a very wide
range, with differences up to $10^{12}$ or even larger. This makes
the constraint matrix very ill-conditioned and the solutions
returned by solvers are often inaccurate and may contain errors.
Also, the presence of  big-\emph{M} terms results in weak bounds
thus leading to very large search trees. To tackle
these problems a number of different approaches were recently
proposed. For a comprehensive introduction to these related works,
we refer the reader to \cite{DA10,DAMaSa10,KeOlRa10}.

In this paper, we propose an alternative pure 0-1 Linear
Programming formulation for the WND, whose defining inequalities
are linear constraints in the assignment variables $x_{tb}$. Such
inequalities are thus valid for all the formulations that are
derived from the previously introduced MILPs and can be included to
strengthen them.

Let now $\tilde{x} \in \{0,1\}^{T\times B}$ be a  server
assignment and let $\Sigma$ denote the set of all the SIR
inequalities $SIR(t,b)$ and the lower and upper bounds
constraints $0 \leq p_b \leq P_{\max}$ on power emissions. With $\tilde{x}$ we associate the subsystem $I(\tilde{x})$ of
SIR inequalities (\ref{eq:firstSIRinequality}) whose corresponding
variables $\tilde{x}_{tb}$ are activated, i.e:
$$
I(\tilde{x}) = \{SIR(t,b) \in \Sigma: \tilde{x}_{tb} = 1\}
$$
It is easy to check  if $I(\tilde{x})$, extended with lower and upper bounds
on the variables $p_b$, is feasible. If this is the case, all of the
assigned testpoints can actually be served by the network, and we
say that $x$ is a {\em feasible server assignment}.

At this point, we can restate the  WND as the problem of finding a
feasible server assignment that maximizes the revenue function. To
this aim, a simple characterization of all the feasible server
assignments goes as follows. Denote by \emph{IS}  the set of subsystems $I(x)$ such that $x$ \emph{is not feasible}.
Then $\tilde{x} \in \{0,1\}^{T\times B}$ is a
feasible server assignment if and only if $\tilde{x}$ satisfies
the following system of linear inequalities:
\begin{equation}\label{eq:WND-IS}
\sum_{(t,b)\in I} \tilde{x}_{tb} \leq |I|-1
\hspace{1.0cm} \forall \hspace{0.1cm} I \in IS
\end{equation}

\noindent
The above system is in general very large and the inequalities must be
generated dynamically. Unfortunately, the separation of violated
inequalities (\ref{eq:WND-IS}) is hard, both theoretically and in practice \cite{AmKa95,Pfetsch_PhD}. Moreover,  it  may entail some of
the numerical difficulties associate with the MILP formulations for
the WND. Still, a relevant subset of these  inequalities
can be separated more effectively, as we describe next.

To this end, we proceed in a similar way to \cite{MaMaSa09}.
Namely, we generate a new system $\Sigma'$ obtained from
$\Sigma$ by substituting each (\ref{eq:firstSIRinequality}) with
the system:
\begin{equation}
\label{eq:binarySIR} \frac{a_{t \beta} \cdot
p_{\beta}}{a_{t b} \cdot p_b} \geq \delta
\hspace{1.0cm} \forall \hspace{0.1cm} b\in B
\setminus\{\beta\}
\end{equation}

\noindent where, to simplify the notation, we assume that $B$ also
contains the noise $\mu$ as a fictitious transmitter with fixed power emission. It is not difficult to see that $\Sigma'$ is a
relaxation of $\Sigma$ and every infeasible subsystem of $\Sigma'$ corresponds to an infeasible subsystem of $\Sigma$.
Basically, this relaxation corresponds to considering a receiver
as served if the power emission of its server suffices to contrast
each interferer individually and the thermal noise. Or,
alternatively, if its {\em best server} is ``stronger'' than its
{\em strongest interferer}. In \cite{MaMaSa09} the authors show that,
in most cases of practical interest, this is indeed a good
approximation of the original SIR constraint.

By assuming $p_b\in [\epsilon, P_b]$, with $\epsilon
> 0$ very small, and by taking the logarithm\footnote{This corresponds to rewriting all quantities in dB format.} of both left and
right hand side multiplied by 10, the system $\Sigma'$  can be
rewritten as:
\begin{equation}\label{eq:differenceSIR}
\label{eq:dBSIR} q_b - q_\beta \leq w_{\beta b}^t \qquad t\in T,
\beta \in B,  b\in B \setminus\{\beta\}
\end{equation}

\noindent
where $q_b = 10\log_{10} p_b$ for all $b\in B$ and $
w^t_{\beta b} = \lceil 10(\log_{10}a_{t \beta} - \log_{10}
a_{t b} -  \log_{10} \delta)\rceil$, extended with the
lower and upper bounds $10\log_{10}\epsilon \leq q_b \leq 10
\log_{10} P_{b}$, for all $b\in B$.

In this way, the system $\Sigma'$ is transformed into a \emph{system
of difference inequalities} (lower and upper bounds
can be easily represented in this form as well), where each
constraint (\ref{eq:differenceSIR}) is associated with a server
$\beta$ and a receiver $t$ and thus with an assignment variable
$x_{t\beta}$.

Now, given a generic system of difference constraints $\Sigma^d$:
\begin{equation}\label{eq:dual-shortest-path}
\begin{array}{lll}
(i) \hskip .5cm & t_v - t_u  \leq l_{uv}, & (u,v) \in A
\\[6pt]
\end{array}
\end{equation}

\noindent where $t\in \reals^A$ and $l \in Z^A $, we can consider
the associated weighted directed graph $G = (V,A)$, with weight
function $l$. Then, it is well known that every infeasible
subsystem of (\ref{eq:dual-shortest-path}) contains (the
constraints corresponding to) the arcs of a negative directed
cycle of $G$ \cite{NeWo88}. Also, denoting by $x\in \{0,1\}^A$ the incidence
vector of (the arcs corresponding to) a feasible subsystem of
$\Sigma^d$, then $x$ is the set of solutions to:
\begin{equation}\label{eq:form-circuit}
\begin{array}{lll}
(i) \hskip .5cm & \sum_{uv\in C} x_{uv} \leq |C| - 1, & C \in
{\cal C^-}
\\[6pt]
 & x\in \{0,1\}^A & \\
\end{array}
\end{equation}

\noindent where ${\cal C^-}$ is the set of negative directed
cycles of $G$.

In \cite{DaMaSa11} we develop an exact approach to the separation
of violated inequalities (\ref{eq:form-circuit}.i).  The resulting
algorithm can be used to separate the violated inequalities
associated with the system (\ref{eq:differenceSIR}) (including
upper and lower bounds on the $q$ variables expressed as
difference inequalities) which correspond to negative directed
cycles in the associated directed graph. One of these cycles $C$
corresponds to a subset of constraints of (\ref{eq:differenceSIR})
associated with the pairs $I_C = \{(\beta_1,t_1), \dots, (\beta_m,
t_m)\}\subseteq B\times T$ (plus possibly some lower and upper
bound constraints).

One can show that $\beta_1 \neq \beta_2 \neq \dots \neq \beta_m$
and $t_1 \neq t_2 \neq \dots \neq t_m$ and the valid constraint:
\begin{equation}\label{eq:WND-cycle}
\sum_{(t,b)\in I_C} x_{tb} \leq |I_C|-1 \hspace{0.5cm}
\end{equation}

\noindent may be added to the formulation. In our preliminary
results, however, we limit to consider  cycle inequalities with
$|C|=2$ and separate them by enumeration.

\section{Preliminary Computational Results}

We test the performance of our new approach to WND on a set of 15
realistic instances, developed with the Technical Strategy \&
Innovations Unit of British Telecom Italia (BT Italia SpA). All
the instances refer  to a \emph{Fixed WiMAX Network}
\cite{AnGhMu07}, deployable in an urban residential area and
consider various scenarios with up to $|T| = 529$ receivers, $|B|
= 36$ transmitters, $|F| = 3$ frequencies, $|H| = 4$ burst
profiles (see Table \ref{tab:comparisons}). We remark that the
experiments refer to a formulation that extends the basic one
considered in Section \ref{sec:intro}, by including frequency
channels and modulation schemes as additional decision variables.
Such formulation is denoted by \emph{(BM)} and captures specific
features of so-called \emph{Next Generation Networks} like WiMAX \cite{AnGhMu07}.
For a detailed description of (BM), we refer the reader to
\cite{DAMa09}.

For each instance, we present preliminary computational results obtained by solving the big-M formulation (BM) and its corresponding Power-Indexed formulation (PI) \cite{DAMaSa10}.
We consider (BM) and (PI) formulations with and without the valid inequalities (\ref{eq:form-circuit}) obtained for $|C| = 2$. Formulations strengthened through (\ref{eq:form-circuit}) are distinguished by adding \emph{S-}, i.e. \emph{(S-BM)} and \emph{(S-PI)}.

Experiments are run by imposing a time limit of 1 hour and by using a machine with a 1.80 GHz Intel Core 2 Duo processor and 2 GB of RAM.
Table \ref{tab:comparisons} reports the performance of the four considered formulations over the set of WiMAX instances. We solve (BM) and (S-BM) by direct application of IBM ILOG Cplex 11.1 and we report i) the upper bound UB$_0$ obtained at node 0 of the branch-and-bound tree, ii) the value $|T^*|$ of the best solution found within the time limit and iii) the final integrality gap \emph{gap\%}. The presence of two values in some lines of the column $|T^*|$ of (BM) indicates that the coverage plans returned by Cplex contain errors and some receivers are actually not covered. We instead solve (PI) and (S-PI) by the incremental algorithm WPLAN described in \cite{DAMaSa10} and we report i) the upper bound UB$_0$ obtained at node 0 when considering the basic set of power levels, and ii) the value $|T^*|$ of the best solution found by WPLAN within the time limit.

By adding the new valid inequalities (\ref{eq:form-circuit}) for
$|C|=2$, in most cases stronger bounds are obtained at node 0 and
smaller integrality gaps are reached within the time limit. In
particular, the benefits are particularly evident in the case of
the big-M formulation: in three cases, namely I6, I12, I15, the
value of the best solution is increased, even eliminating coverage
errors (I6, I15).

\begin{table}[tbp]
\caption{Comparisons between (BM) and (PI) with and without valid inequalities (\ref{eq:form-circuit})}
\label{tab:comparisons} \small
\begin{center}
\begin{tabular}{|ccccc||ccc|ccc|ccc|ccc|}
  \hline
  & & & & &
  %&
  & (BM) & &
  & (S-BM) & &
  & (PI) & &
  & (S-PI) &
  \\
  \raisebox{1.5ex}{ ID } & \raisebox{1.5ex}{$|\mbox{T}|$} &
  \raisebox{1.5ex}{$|\mbox{B}|$} &
  \raisebox{1.5ex}{$|\mbox{F}|$} &
  \raisebox{1.5ex}{$|\mbox{H}|$} &
  UB$_0$ & $|T^*|$ & gap\% &
  UB$_0$ & $|T^*|$ & gap\% &
  UB$_0$ && $|T^*|$ &
  UB$_0$ && $|T^*|$
%  $\mbox{ }$ WPLAN \cite{DAMaSa10} $\mbox{ }$ &
%  (GA-PFMAP) $\mbox{ }$
  \\
  \hline
  \hline
  I1 & 100 & 12 & 1 & 1
  & 98.36 & 66 (70) & 29.43
  & $\mbox{ }$ 96.78 & 66 (70) & 27.68 $\mbox{ }$
  & 90.77 && 75
  & $\mbox{ }$ 90.23 && 75
  \\
  I2 & 169 & 12 & 1 & 1
  & 165.47 & 97 & 59.81
  & $\mbox{ }$ 163.15 & 97 & 57.39 $\mbox{ }$
  & 153.12 && 101
  & $\mbox{ }$ 152.37 && 101
  \\
  I3 & 196 & 12 & 1 & 1
  & 193.61 & $\mbox{ }$ 102 (105) $\mbox{ }$ & 77.87
  & $\mbox{ }$ 192.02 & $\mbox{ }$ 102 (105) $\mbox{ }$ & 75.11 $\mbox{ }$
  & 179.35 && 108
  & $\mbox{ }$ 177.92 && 108
  \\
  I4 & 225 & 12 & 1 & 1
  & 219.76 & 92 & 81.13
  & $\mbox{ }$ 218.36 & 92 & 79.36 $\mbox{ }$
  & 202.44 && 92
  & $\mbox{ }$ 201.54 && 92
  \\
  I5 & 289 & 12 & 1 & 1
  & 287.20 & 76 (77) & 195.44
  & $\mbox{ }$ 287.20 & 76 (77) & 194.92 $\mbox{ }$
  & 274.62 && 85
  & $\mbox{ }$ 274.13 && 85
  \\
  I6 & 361 & 12 & 1 & 1
  & 352.01 & 126 (132) & 154.87
  & $\mbox{ }$ 350.43 & 140 & 138.76 $\mbox{ }$
  & 337.22 && 156
  & $\mbox{ }$ 336.46 && 156
  \\
  I7 & 400 & 18 & 1 & 1
  & 397.21 & 166 & 132.01
  & $\mbox{ }$ 396.79 & 166 & 131.32 $\mbox{ }$
  & 386.07 && 184
  & $\mbox{ }$ 384.95 && 184
  \\
  I8 & 400 & 18 & 3 & 4
  & 400.00 & 356 & 12.36
  & $\mbox{ }$ 400.00 & 356 & 12.36 $\mbox{ }$
  & 396.53 && 372
  & $\mbox{ }$ 395.80 && 372
  \\
  I9 & 441 & 18 & 3 & 4
  & 441.00 & 266 (270) & 63.33
  & $\mbox{ }$ 441.00 & 266 (270) & 63.33 $\mbox{ }$
  & 438.28 && 295
  & $\mbox{ }$ 437.52 && 295
  \\
  I10 & 484 & 27 & 3 & 4
  & 484.00 & 120 (122) & 296.72
  & $\mbox{ }$ 484.00 & 120 (122) & 296.72 $\mbox{ }$
  & 479.10 && 242
  & $\mbox{ }$ 478.68 && 242
  \\
  I11 & 529 & 27 & 3 & 4
  & 529.00 & 77 & 587
  & $\mbox{ }$ 529.00 & 77 & 587 $\mbox{ }$
  & 523.15 && 168
  & $\mbox{ }$ 521.76 && 168
  \\
  I12 & 400 & 36 & 1 & 4
  & 398.04 & 72 (74) & 287.30
  & $\mbox{ }$ 396.93 & 77 (78) & 264.85 $\mbox{ }$
  & 389.61 && 102
  & $\mbox{ }$ 389.14 && 102
  \\
  I13 & 441 & 36 & 1 & 4
  & 433.21 & 184 & 131.03
  & $\mbox{ }$ 431.42 & 184 & 129.77 $\mbox{ }$
  & 414.93 && 194
  & $\mbox{ }$ 413.78 && 194
  \\
  I14 & 484 & 36 & 1 & 4
  & 482.78 & 209 & 108.31
  & $\mbox{ }$ 481.66 & 209 & 107.56 $\mbox{ }$
  & 472.44 && 251
  & $\mbox{ }$ 471.58 && 251
  \\
  I15 & 529 & 36 & 1 & 4
  & 517.89 & 98 (105) & 226.44
  & $\mbox{ }$ 516.14 & 114 & 198.57 $\mbox{ }$
  & 503.32 && 232
  & $\mbox{ }$ 502.67 && 232
  \\
\hline
\end{tabular}
\end{center}
\end{table}

\end{document}